\documentclass[11pt,reqno]{amsart}
\usepackage{amssymb}
\usepackage{cite}
\usepackage{mathrsfs}
\usepackage{amsfonts}
\usepackage{amsmath}
\def\MM{\ensuremath{{\mathscr M}}}
\def\pv#1{\ensuremath{{\bf#1}}}

\def\inv{^{-1}}
\def\p{\varphi}

\def\pinv{{\p \inv}}
\def\J{\mathrel{{\mathscr J}}} 

\def\R{\mathrel{{\mathscr R}}} 
\def\L{\mathrel{{\mathscr L}}} 
\def\H{\mathrel{{\mathscr H}}} 

\def\e<{\leq _{E}}

\def\ov#1{\ensuremath{\overline {#1}}}

\def\pvr#1{\ensuremath{\mathsf{#1}}}

\def\1sk{^{(1)}}

\def\GGM{\mathsf{GGM}}

\def\Rad#1#2{\rho_{\pv {#1}}(#2)}

\def\to{\rightarrow}

\def\malce{\mathbin{\hbox{$\bigcirc$\rlap{\kern-8.75pt\raise0,50pt\hbox{$\mathtt{m}$}}}}}
\def\malced{\mathbin{\hbox{$\bigcirc$\rlap{\kern-7pt\raise0,50pt\hbox{$\mathtt{m}$}}}}}
\def\data{\ifcase\month\or January\or February \or March\or April\or May
\or June\or July\or August\or September\or October\or November \or
December\fi\space\number\day, \number\year}

\def\Thmname{Theorem}
\def\Propname{Proposition}
\def\Lemmaname{Lemma}
\def\Definitionname{Definition}
\newtheorem{Thm}{\Thmname}[section]
\newtheorem{Prop}[Thm]{\Propname}
\newtheorem{Lemma}[Thm]{\Lemmaname}
\newtheorem{Def}[Thm]{\Definitionname}

\newtheorem{Cor}[Thm]{Corollary}

\numberwithin{equation}{section}

\title[A structural approach to locality]{A structural approach to the
  locality of pseudovarieties of the form $\pv {LH}\malce \pv V$}

\author   {B.~Steinberg} \address{School of Mathematics and
  Statistics, Carleton
University, Ottawa, Ontario K1S 5B6, Canada}
\email{bsteinbg@math.carleton.ca}
\thanks{The author gratefully acknowledge the support of NSERC}

\dedicatory{Dedicated to the memory of Bret Tilson}

\date{December 16, 2006}


\begin{document}

\begin{abstract}
We show that if \pv H is a Fitting pseudovariety of groups and \pv V
is a local pseudovariety of monoids, then $\pv {LH}\malced \pv V$ is
local if either $\pv V$ contains the six element Brandt monoid, or
\pv H is a non-trivial pseudovariety of groups closed under
extension.
\end{abstract}

\maketitle

\section{Introduction}
Since the seminal work of Tilson \cite{Tilson}, as well as further
work of Tilson and Rhodes \cite{Kernel}, Rhodes and Weil
\cite{RWI,RWII} and work of Pin, Straubing and Therien \cite{PST},
it has become clear that finite categories play a crucial role in
finite monoid and formal language theory.  Tilson \cite{Tilson}
defined a pseudovariety of monoids \pv V to be \emph{local} if the
pseudovariety of categories generated by the elements of \pv V,
viewed as one-object categories, consists precisely of those
categories whose local monoids  at each
object belong to \pv V.  Many important pseudovarieties of monoids
are local, but not all; see \cite{Tilson} for more background.

Let $\pv H$ be a pseudovariety of groups.   A pseudovariety of
groups $\pv H$ is said to be \emph{Fitting}~\cite{Fitting} if
whenever $G$ is a finite group and $H,K\lhd G$ are normal subgroups
belonging to $\pv H$, then $HK\in \pv H$. For instance,
pseudovarieties of groups closed under extension are Fitting; so is
the pseudovariety of nilpotent groups. Notice that if \pv H is
Fitting and $G$ is a finite group, then there is a largest normal
subgroup $\Rad H G$ of $G$ belonging to \pv H, namely the product of
all normal subgroups $N\lhd G$ with $N\in \pv H$. This subgroup is
called the \emph{$\pv H$-radical} of $G$.  For example, the
$p$-group radical of a finite group $G$ is just the intersection of
all its $p$-Sylow subgroups; the nilpotent radical of $G$ is the
product of the $p$-group radicals taken over all prime divisors $p$
of the order of $G$.  For more on Fitting pseudovarieties of groups
and their relationship with semigroup theory, see~\cite{Fitting}.

We denote by $\pv {LH}$ the pseudovariety of all semigroups $S$ whose
\emph{local monoids} $eSe$, with $e$ an idempotent, belong
to $\pv H$.  If \pv V is a pseudovariety of semigroups and \pv W
is a pseudovariety of monoids, then the \emph{Mal'cev product}
$\pv V\malce \pv W$ is the pseudovariety generated by all
monoids $M$ with a homomorphism $\p:M\to N\in \pv W$ such that
$e\pinv \in \pv V$ for each idempotent $e\in N$.  It in fact
consists precisely of all homomorphic images of such monoids $M$.

In this paper, we shall be interested in Mal'cev products \mbox{$\pv
{LH}\malce \pv V$} where \pv H is a Fitting pseudovariety of groups.
Let $\pv {G}_p$ be the pseudovariety of $p$-groups, where $p$ is
prime, and let $\pv {G_{sol}}$ be the pseudovariety of solvable
groups:  both are Fitting. The work of Weil~\cite{Weil} describes
the languages in $\pv L\pv {G}_p\malce \pv V$ and $\pv L\pv
{G_{sol}}\malce \pv V$ in terms of the languages of \pv V using
marked products with modular counters. The trivial pseudovariety of
groups \pv I is also Fitting. It is shown in \cite{PST} that the
languages of $\pv {LI}\malce \pv V$ can be described in terms of
unambiguous marked products of \pv V-languages.

Our main theorem is the following.

\begin{Thm}\label{mainthm}
Let \pv H be a pseudovariety of groups and \pv V a local pseudovariety
of monoids.  Then $\pv{LH}\malce \pv V$ is local under either of the
following circumstances:
\begin{enumerate}
\item  \pv H is a Fitting pseudovariety of groups and \pv V contains
  the 6 element Brandt monoid $B_2^1$;
\item \pv H is non-trivial and closed under extension.
\end{enumerate}
\end{Thm}

The proof follows the scheme of \cite[Section 6]{assertion}, where
locality is established for a one-sided version of $\pv {LG}\malce
\pv V$ (here \pv G is the pseudovariety of all finite groups).
We also use ideas from
\cite{Kernel,RWII,mystiffler,JT,therien2side}.

Let \pv{DS} be the pseudovariety of monoids whose regular
$\J$-classes are subsemigroups.  If \pv H is a pseudovariety of
groups, then $\ov {\pv H}$ denotes the pseudovariety of monoids
whose subgroups belong to \pv H. Consider the pseudovariety $\pv
{DS}\cap \ov {\pv H}$.  Independent work of Putcha and
Sch\"utzenberger \cite{Putcha,Schutzambi} (see also \cite{Almeida})
shows that $$\pv {DS}\cap \ov {\pv H} = \pv {LH}\malce \pv {Sl}.$$
Since \pv {Sl} is local by a theorem of Simon
\cite{Simon,Eilenberg}, Theorem \ref{mainthm}(2) shows  $\pv
{DS}\cap \ov {\pv H}$ is local when \pv H is extension-closed and
non-trivial, a result first proved in \cite{JT}.  The case \pv H is
trivial is the pseudovariety \pv{DA}, which was shown to be local by
Almeida \cite{AlmeidaDA}.  This was first announced by Therien
\cite{therien2side}, but this paper contains an unproven lemma,
which was left as an exercise and which is generally believed to be
an open problem. If this lemma were true, our techniques would show
that Theorem \ref{mainthm} applies in the case \pv H is trivial
without the assumptions of (1).

If \pv {DO} is the pseudovariety of monoids from \pv {DS} whose
idempotents in each $\J$-class form a subsemigroup, then each
subpseudovariety of $\pv {DO}$  whose associated variety of languages
is closed under
unambiguous product was shown to be local in \cite{AlmeidaDA} using
profinite techniques.  The results of \cite{PST} imply that any such
pseudovariety is of the form $\pv {LI}\malce \pv V$ where $\pv V$ is a
pseudovariety of semilattices of groups; see \cite[Corollary
  3.3]{AlmeidaDA}.  Such pseudovarieties \pv V are
local by work of Jones and Szendrei \cite{JS}, so Theorem
\ref{mainthm} is true in these cases even when \pv H is trivial
(without the $B_2^1$ assumption). We leave it as an open question as
to whether the assumption that \pv H is non-trivial is needed in
Theorem \ref{mainthm}(2). More generally, does Theorem \ref{mainthm}
(2) hold for any Fitting pseudovariety of groups (i.e.\ is the
hypothesis on \pv V in Theorem \ref{mainthm}(1) necessary)?

Our paper is organized as follows.  We first consider \pv
{LH}-morphisms of categories.  We show that any category admitting an
\pv {LH}-morphism to an element of $\pv {gV}$ belongs to $\pv g(\pv
    {LH}\malce \pv V)$ under the assumptions of Theorem \ref{mainthm}.
    Afterwards, we show that if
    \pv V is local, then any category that is locally in $\pv
    {LH}\malce \pv V$ admits an \pv {LH}-morphism to an element of
    \pv{gV}.  This is our technically most difficult theorem and
    requires the use of radical congruences from
    \cite{Arbib,folleyR,folleyT,HW,Fitting} and generalized group mapping
    monoids~\cite{Arbib,qtheor}.  In order to keep this paper
    as short as possible, we prove only what is necessary.  In
    particular, we don't fully develop the theory of $\pv
    {LH}$-morphisms of categories, where most results from the
    monoid setting can be shown to hold.

In this paper all semigroups, monoids, categories and semigroupoids
are finite.

\section{\pv {LH}-morphisms of categories}
We begin by defining \pv{LH}-morphisms.  The reader is referred to
\cite{Tilson,MacLane} for the notion of a quotient morphism of
categories.   If $C$ is a
category, $C(c,c')$ denotes the hom set of arrows
from $c$ to $c'$.  We use $C_c$ for the local monoid $C(c,c)$.  The
set of arrows of $C$ is denoted $\pvr{Arr}(C)$; the set of objects of
$C$ is denoted $\pvr{Obj}(C)$.   If
$m\in \pvr{Arr}(C)$, then $m\iota$ is the initial vertex of $m$ and
$m\tau$ is the terminal vertex.  We write $m:m\iota\to m\tau$.

\begin{Def}[$\pv {LH}$-morphism]
A morphism of categories $\p:C\to D$ is called an \emph{$\pv
  {LH}$-morphism} if it is a quotient morphism and if, for each
  idempotent $e\in \pvr {Arr} (D)$, $e\pinv\in \pv {LH}$.
\end{Def}

Since the congruence class of an idempotent is contained in a
local monoid, the notion of an $\pv {LH}$-morphism is local in the
following sense.

\begin{Prop}\label{forgotten}
Let $\p:C\to D$ be a quotient morphism of categories.  Then $\p$
is an \pv{LH}-morphism if and only if, for all $c\in
\pvr{Obj}(C)$, $\p|_{C_c}:C_c\to D_c$ is an \pv {LH}-morphism (of monoids).
\end{Prop}

If \pv V is a pseudovariety of monoids,
the  pseudovariety $\pv {LH}\malce \pv V$ consists of all quotients of
monoids admitting an \pv {LH}-morphism to an element of \pv V.
If \pv H is a Fitting pseudovariety of groups, then $\pv {LH}\malce
\pv V$ consists of all monoids admitting an \pv {LH}-morphism to a
monoid in \pv V; see Section \ref{awfultech}.

The following proposition generalizes a result from the monoid case.
The converse
is true if \pv H is extension-closed, although we do not prove it as
we shall not use it.

\begin{Prop}\label{comp}
Suppose $\p:C\to D$ is an $\pv{LH}$-morphism and $\p = \psi\gamma$ where
$\psi:C\to C'$ and $\gamma:C'\to D$ are quotient morphisms.  Then $\psi$ and
$\gamma$ are \pv {LH}-morphisms.
\end{Prop}
\begin{proof}
By Proposition \ref{forgotten}, it suffices to show that for each
$c\in \pvr{Obj}(C)$, $\psi|_c$ and $\gamma|_c$ are \pv
{LH}-morphisms.  Since $\p_c$ is an \pv {LH}-morphism by Proposition
\ref{forgotten}, we have that  $\psi|_c$ and $\gamma|_c$ are \pv
{LH}-morphisms by the well-known analogue of Proposition \ref{comp}
for monoids \cite{Kernel,RWII}.
\end{proof}

A key idea in this paper is to translate things from categories to
monoids via the consolidation operator.
If $C$ is a category, we denote by $C^{cd}$ the \emph{consolidation} of $C$.
It is the monoid $C\cup 0\cup 1$ obtained by adjoining a zero $0$
and an identity $1$ to $C$ and making all previously undefined
products zero.

It was observed in \cite{mystiffler} that if $\p:C\to D$ is a
quotient morphism of categories, then there is an induced monoid
morphism $\p_{cd}:C^{cd}\to D^{cd}$ given by $\p$ on
$\pvr{Arr}(C)$ and  sending $0$ to $0$ and $1$ to $1$.
The next proposition shows that the notion of $\pv {LH}$-morphism carries over
nicely to the consolidation.

\begin{Prop}\label{LHtocd}
 Let $\p:C\to D$ be a quotient morphism.  Then $\p$ is an \pv
  {LH}-morphism if and only if $\p_{cd}:C^{cd}\to D^{cd}$ is an \pv
  {LH}-morphism.
\end{Prop}
\begin{proof}
First note that the idempotents of $D^{cd}$ are the idempotents of $D$
along with $0$ and $1$.  Moreover, if $e\in \pvr {Arr}(D)$ is an
idempotent, then $e\pinv = e\p_{cd}\inv$.
Since $1\p_{cd}\inv =1$, $0\p_{cd}\inv =0$  it now is clear that $\p_{cd}$
is an \pv {LH}-morphism if and only if $e\pinv \in \pv {LH}$ for each
idempotent $e\in \pvr {Arr}(D)$ if and only if $\p$ is an \pv {LH}-morphism.
\end{proof}

To prove Theorem \ref{mainthm}, we need the following
technical result, whose proof we defer to Section \ref{awfultech}.

\begin{Thm}\label{supertech}
Let $\pv V$ be a pseudvariety of monoids and \pv {H} be a Fitting
pseudovariety of groups.  Let $C$ be a category.
Then $C\in \ell (\pv {LH}\malce \pv V)$ if and only if $C$ admits an
\pv {LH}-morphism to a memeber of $\ell \pv V$.
\end{Thm}

\subsection{Proof of Theorem \ref{mainthm}(1)}
Suppose that $C\in \ell (\pv {LH}\malce \pv V)$.  Then there is,
by Theorem \ref{supertech}, an $\pv{LH}$-morphism $\p:C\to D$ with
$D\in \ell\pv V= \pv{gV}$.  Hence, by Proposition \ref{LHtocd},
$\p_{cd}:C^{cd}\to D^{cd}$ is an $\pv {LH}$-morphism.  But since
$D\in \pv {gV}$, we have $D^{cd}\in \pv V$, as $\pv V$ contains $B_2^1$
\cite{Tilson}.  Thus $C^{cd}\in \pv {LH}\malce \pv V$, whence
$C\in \pv g(\pv {LH}\malce \pv V)$, since $C$ divides
$C^{cd}$.\qed

\subsection{Proof of Theorem \ref{mainthm}(2)}
We now recall the notion of a maximal proper quotient (MPQ)
\cite{mystiffler}, generalizing Rhodes' notion of a maximal proper
surmorphism (MPS) \cite{Rhodesmps,RWI,Arbib}.  This notion will play
a key role in the proof of Theorem \ref{mainthm}(2).

\begin{Def}[Maximal proper quotient] A quotient morphism $\p:C\to D$
  of categories is called a \emph{maximal proper quotient} (MPQ), if
  the associated congruence $(\p)$ is minimal amongst non-trivial
  congruences on $C$.
\end{Def}

A monoid morphism that is an MPQ is called a \emph{maximal proper
  surmorphism} or MPS.

Clearly any quotient morphism $\p:C\to D$ of categories factors
$\p = \p_1\cdots \p_n$ with each $\p_i$ an MPQ.  This combined with
Proposition \ref{comp} proves.

\begin{Prop}\label{mpsreduction}
Let \pv H be a pseudovariety of groups.  Let $\p:C\to D$ be a \pv
{LH}-morphism.  Then $\p = \p_1\cdots \p_n$ where each $\p_i$ is both
an MPQ and a \pv {LH}-morphism.
\end{Prop}

To prove Theorem \ref{mainthm}(2), we shall essentially use an
induction argument on the length of the above decomposition. In
\cite{mystiffler}, we made an observation that allows us reduce many
results about MPQs to MPSs.  This is done via the consolidation
operator.   The following is \cite[Proposition 2.1]{mystiffler}.

\begin{Prop}\label{toconsolidate}
  Let $\p:C\to D$ be a quotient morphism.  Then $\p$ is a maximal
  proper quotient if and only if $\p_{cd}:C^{cd}\to D^{cd}$ is a
  maximal proper surmorphism.
\end{Prop}

In order to use some decomposition results from semigroup theory,
we shall need to consider the kernel category of a
quotient morphism of categories \cite{JP,myCats}.
 Let $\p:C\to D$ be
a quotient morphism; we view $\p$ as being the identity on
$\pvr{Obj}(C)$.   First we define a category $W_{\p}$.  Its object
set consists of all pairs $(n_L,n_R)$ of arrows of $D$ such that
$n_L\tau = n_R\iota$ (that is, $n_Ln_R$ is defined).  The arrows are
of the form
$$(n_L,m,n_R):(n_L,m\p n_R)\to (n_Lm\p, n_R),$$ where $n_L\tau= m\iota$
and $m\tau = n_R\iota$ (so $n_Lm\p n_R$ is defined).  Multiplication
is given by
$$(n_L,m,m'\p n_R)(n_Lm\p,m',n_R) = (n_L,mm',n_R).$$
The identity at $(n_L,n_R)$ is $(n_L,1_{n_L\tau},n_R)$.  The kernel category
\cite{JP,myCats}
$K_{\p}$ is the quotient of $W_{\p}$ by the congruence that identifies
two coterminal
arrows $(n_L,m,n_R)$ and $(n_L,m',n_R)$ if and only if, for all
$m_L\in n_L\pinv$, \mbox{$m_R\in n_R\pinv$},
\begin{equation}\label{ident}m_Lmm_R=m_Lm'm_R
\end{equation}
  When $C$
and $D$ are monoids, this is the kernel category of \cite{Kernel}.

Our next theorem is a reformulation of a result that was first
announced in \cite{JP} (see also \cite{therien2side}), but was first
correctly proved by the author \cite{myCats} (see \cite{Cats2} for
the one-sided analogue).  This result plays a key role in \cite{JT}.
Recall \cite{Kernel} that $\pv V\mathrel{\ast\ast} \pv W$ denotes
the two-sided semidirect product of pseudovarieties.

\begin{Thm}\label{globalthm}
Let \pv V and \pv W be pseudovarieties of monoids and let $\p:C\to
D$ be a quotient morphism with $D\in \pv {gW}$ and $K_{\p}\in \pv
{gV}$. Then $C\in \pv g(\pv V\mathrel{\ast\ast} \pv W)$.
\end{Thm}

To use the powerful decomposition results of Rhodes \emph{et al.}
\cite{Kernel,RWII}, we need to relate the $K_{\p}$ and $K_{\p_{cd}}$.

\begin{Prop}\label{kernelin}
Let $\p:C\to D$ be a quotient morphism of categories.  Then $K_{\p}$
is the full subcategory of $K_{\p_{cd}}$ with objects $(n_L,n_R)$ with
$n_L, n_R\in \pvr{Arr}(D)$, $n_L\tau = n_R\iota$.  Hence, $K_{\p}$
belongs to any
pseudovariety of categories containing $K_{\p_{cd}}$.
\end{Prop}
\begin{proof}
Clearly $W_{\p}$ embeds in $W_{\p_{cd}}$ as a subsemigroupoid,
however the local identities are different.  The only arrows of
$W_{\p_{cd}}$ between objects of $W_{\p}$ that do not belong to
$W_{\p}$ are the arrows of the form $(n_L,1,n_R)$, where $n_L\tau =
n_R\iota$. But \eqref{ident} shows that if $(n_L,n_R)$ is in
$\pvr{Obj}(K_{\p})$, then  $(n_L,1,n_R)$ is identified with
$(n_L,1_{n_L\tau},n_R)$ in $K_{\p_{cd}}$.   Also it is clear from
\eqref{ident} that two arrows of $W_{\p_{cd}}$ become identitied in
$K_{\p_{cd}}$ if and only if they are identified in $K_{\p}$.
Putting this all together we see $K_{\p}$ is a full subcategory of
$K_{\p_{cd}}$.
\end{proof}

We now pull out a big hammer, namely the following deep result of
Rhodes \emph{et al.} \cite{Kernel,RWII}.

\begin{Thm}\label{mpsresult}
Let \pv H be a pseudovariety of groups.  Let $\p:M\to N$ be a maximal
proper surmorphism of
monoids and an \pv{LH}-morphism.  Then $K_{\p}\in \pv \ell\pv H$.
\end{Thm}

Putting everything together we obtain:

\begin{Cor}\label{passtocat}
Let \pv H be a pseudovariety of groups and let $\p:C\to D$ be a
maximal proper quotient
and an \pv{LH}-morphism of categories.  Then $K_{\p}\in \ell \pv H$.
\end{Cor}
\begin{proof}
Since $\p$ is MPQ, $\p_{cd}$ is MPS by Proposition
\ref{toconsolidate}.  Since $\p$ is an $\pv {LH}$-morphism, so is
$\p_{cd}$ by Proposition \ref{LHtocd}.  By Theorem \ref{mpsresult},
$K_{\p_{cd}}\in \ell\pv H$.  Hence by Proposition \ref{kernelin},
$K_{\p}\in \ell\pv H$.
\end{proof}

We proceed to prove Theorem \ref{mainthm}(2) assuming Theorem
\ref{supertech}.
\subsubsection{Proof of Theorem \ref{mainthm}(2)}
Fix a non-trivial extension-closed pseudovariety of group \pv H.
We first use the following theorem \cite{Kernel,RWII} (see also
\cite{Weil}) whose difficult direction relies on Theorem
\ref{mpsresult} and the fact that $\ell \pv H=\pv {gH}$
\cite{Tilson} if $\pv H$ is non-trivial.

\begin{Thm}\label{LHtheorem}
Let \pv V be a pseudovariety of monoids and \pv H be a non-trivial
extension closed pseudovariety of groups.  Then $\pv {LH}\malce \pv
V$ is the smallest pseudovariety \pv W containing \pv V such that
$\pv H\mathrel{\ast\ast} \pv W = \pv W$.
\end{Thm}

Suppose now that $\pv V$ is a local pseudovariety of monoids.  Let
$C\in \ell (\pv {LH}\malce \pv V)$.  By Theorem \ref{supertech}, $C$
admits an \pv {LH}-morphism $\p:C\to D$ with $D\in \ell \pv V= \pv
{gV}$. Hence by Lemma \ref{mpsreduction} there is a sequence
$\p_1,\ldots,\p_n$ of $\pv {LH}$-morphisms that are MPQ such that
$C$ is the domain $\p_1$, the codomain of $\p_n$ belongs to $\pv
{gV}$ and $\p_1\p_2\cdots \p_n$ is defined.  We proceed by
induction on $n$.  If $n=0$, then $\p$ is the identity map and so
$C\in \pv {gV}$ and we are done.  Else, suppose $\p_1:C\to C'$.
Then $C'\in \ell (\pv {LH}\malce \pv V)$ and by considering
$\p_2,\ldots, \p_n$ and the induction hypothesis, $C'\in \pv g(\pv
{LH}\malce \pv V)$.

By Corollary \ref{passtocat}, $K_{\p_1}\in \ell \pv H = \pv {gH}$, the
equality following from \cite{Tilson}.  So by Theorem \ref{globalthm}
and Theorem \ref{LHtheorem}
$$C  \in \pv g (\pv H\mathrel{\ast\ast} (\pv {LH}\malce \pv V))
      = \pv g (\pv {LH}\malce \pv V).$$  This establishes Theorem
      \ref{mainthm}~(2).\qed

\section{Radical Congruences}\label{awfultech}
The goal of this section is to prove Theorem \ref{supertech}.   In
particular, we show how certain radical congruences on $C^{cd}$ restrict
nicely to local monoids $C_c$.
 We
shall need to use the description of the maximal congruence
associated to an $\pv{LH}$-morphism of monoids for Fitting
pseudovarieties of groups.

 A semigroup $S$ is
  called \emph{generalized group
  mapping}~\cite{Arbib,qtheor} (GGM) if it has a $0$-minimal ideal $I$ on
  which it acts
  faithfully on both the left and right.  
  GGM semigroups play a
  crucial role in finite
  semigroup theory, both in relation to complexity
  theory~\cite{Arbib,folleyR,folleyT,qtheor} and to Malcev
  products of the form $\pv {LH}\malce \pv V$~\cite{qtheor,Fitting}.

We recall some notions and results of Krohn and Rhodes.  The reader
is referred to~\cite{Arbib,qtheor} for details.  Fix a finite
semigroup $S$. Let us set some notation. Choose, for each regular
$\J$-class $J$, a fixed maximal subgroup
$G_J$. 
%
%
Let $J$ be a regular $\J$-class of $S$.
  Let
\[F(J) = \{s\in S\mid s\not\geq_{\J} J\}.\]  Then $F(J)\cup J$ and
$F(J)$ are ideals of $S$.  We identify $(F(J)\cup J)/F(J)$ with
$J^0$. Choose an isomorphism of $J^0$ with a Rees matrix semigroup
$\MM^0(G_J,A,B,C)$.  We do not distinguish between $J^0$ and this
Rees matrix semigroup group.   Now let $N\lhd G_J$ be a normal
subgroup.  Then we let $\eta:F(J)\cup J\to J^0$ and
$\psi:\MM^0(G_J,A,B,C)\to \MM^0(G_J/N,A,B,\ov C)$ be the natural
projections, where $\ov C$ is the matrix $C$ reduced modulo $N$.

We define a congruence $\equiv_{(J,G_J,N)}$ as follows \cite{Arbib}.
Let $s,t\in S$ and $x,y\in J$.  Then $xsy,xty\in F(J)\cup J$.  Define,
for $s,t\in S$,
\begin{equation}\label{radcong}
s\equiv_{(J,G_J,N)}t \iff xsy\eta\psi = xty\eta\psi\ \text{for all}\
x,y\in J.
\end{equation}
The quotient $S/\!{\equiv_{(J,G_J,N)}}$ is denoted $\GGM
(J,G_J,N)$~\cite{Arbib}.   It can be shown that $\GGM (J,G_J,N)$  is
GGM, its  definition does not depend on the choice of the Rees
matrix representation of $J^0$ and that all generalized group
mapping images of $S$ are of this form for some regular $\J$-class
$J$ and normal subgroup $N$ of $G_J$; see~\cite[Proposition 8.3.28,
Remark 8.3.29]{Arbib}.

Let us make an observation about the definition of
$\equiv_{(J,G_J,N)}$. Namely, let $x,y\in J$.  Let $x',y'$ be inverses
of $x,y$ respectively.
Let $s,t\in S$.  Then it is clear that \eqref{radcong} holds if
and only if $$x'xsyy'\eta\psi = x'xtyy'\eta\psi.$$
Thus in \eqref{radcong}, we may always
assume that $x$ and $y$ are idempotent.  We shall use this
observation without further remark.

%


Proofs of the following results are implicit in
\cite{Arbib,folleyR,folleyT}, where they are proved for every
Fitting pseudovariety of groups that had been considered in
semigroup theory up until the time they were written;  notice these
works were written before Eilenberg and Sch\"utzenberger introduced
the notion of pseudovarieties to semigroup theory
\cite{Eilenberg,Eilenschutz} and before anybody had considered
Mal'cev products in the theory.  The general case, whose proof is no
different than the case of $p$-groups already considered by Rhodes
and Tilson \cite{folleyR,folleyT}, can be found in Hall and Weil
\cite{HW}; see also~\cite{Fitting}.  We state the results in the
form that we shall use them.

\begin{Thm}[Rhodes, Tilson]\label{Prop:LHmap}
Let \pv H be a Fitting pseudovariety of groups and $S$ a
semigroup. Then $\p:S\to \prod \GGM (J,G_J,\Rad H {G_J})$, the
product of the canonical morphisms (where the product runs over
all regular $\J$-classes), is an \pv {LH}-morphism and the congruence
associated to $\p$ contains the congruence associated to any other \pv
{LH}-morphism.
\end{Thm}

\begin{Thm}[Rhodes, Tilson]\label{Thm:GGMLH}
Let $\pv H$ be a Fitting pseudovariety of groups, \pv V a
pseudovariety of monoids and $M$ a finite monoid.  Then $M$ belongs
to $\pv {LH}\malce \pv V$ if and only if one has
\begin{equation*}\label{eqofthmggmlh}
\GGM (J,G_J,\Rad H {G_J})\in \pv V
\end{equation*}
for all regular $\J$-classes $J$ of $M$.
\end{Thm}

We are primarily interested in the following corollary.

\begin{Cor}
Let \pv H be a Fitting pseudovariety of groups and let \pv V be a
pseudovariety of monoids.  Then $M\in \pv {LH}\malce \pv V$ if and
only if $M$ has a $\pv {LH}$-morphism to an element of \pv V. In
particular if \pv H and \pv V are decidable, then so is $\pv
{LH}\malce \pv V$.
\end{Cor}

The next two results are analogous to \cite[Lemma 6.2]{assertion}.
Let $C$ be a  category and $c\in \pvr{Obj}(C)$.   To
prove Theorem \ref{supertech}, we must understand how the
$\J$-relation on $C^{cd}$ and $C_c$ relate.

\begin{Prop}\label{Jrelpasses}
Let $x,y\in C_c$.  Then $x\J y$ in $C^{cd}$ if and only if $x\J y$
in $C_c$.
\end{Prop}
\begin{proof}
Clearly $x\J y$ in $C_c$ implies $x\J y$ in $C^{cd}$.  For the
converse, observe that if $uxv =y$ with $u,v\in C^{cd}$, then if
$u$ is not $1$, it must start at $y\iota$ and end at
$x\tau$.  Since $c=y\iota = x\tau$, we see $u\in C_c$.  Similarly,
$v=1$ or $v\in C_c$.  Hence $x\J y$ in $C_c$.
\end{proof}

\begin{Prop}\label{maxsubgroups}
Let $x\in C_c$.  Then the $\H$-class of $x$ in $C^{cd}$ and $C_c$
coincide.  In particular, if $e\in C^{cd}$ is an idempotent other than
$0$ and $1$, then $e\in C_c$ for some $c\in \pvr{Obj}(C)$ and the maximal
subgroup at $e$ in $C^{cd}$ and $C_c$ coincide.
\end{Prop}
\begin{proof}
The second statement is an immediate consequence of the first.  For
the first, observe that in a category $\R$-equivalent elements must
start at the same vertex and $\L$-equivalent elements must end at the
same vertex.  Hence $\H$-equivalent elements are coterminal.  Suppose
now that $x\in C_c$. Then its $\H$-class is contained in $C_c$  The
fact that the $\H$-class of $x$ in
$C^{cd}$ and $C_c$ coincide now follows from \cite[Lemma 6.2]{assertion}.
\end{proof}

These two propositions can be summarized as follows:

\begin{Cor}\label{unionofH}
Let $C$ be a category and $c\in \pvr{Obj}(C)$.  Let $J$ be a
$\J$-class of $C^{cd}$.  If $J_c=J\cap C_c$ is not empty, then it  is a
$\J$-class of $C_c$ and
is a union of $\H$-classes of $J$.  Moreover, $J$ is regular if and
only if $J_c$ is regular.
\end{Cor}
\begin{proof}
Only the last statement has not been proved.  Clearly if $J_c$ is
regular, so is $J$.  Suppose $J$ is regular.   Let $x\in J_c$.  Let
$x'$ be an inverse of $x$.  Then $xx'x = x$ implies that $x'\in
C_c\cap J = J_c$.  Thus $J_c$ is regular.
\end{proof}

Let $C$ be a category and fix a maximal subgroup $G_J$ for
each regular $\J$-class $J$ of $C^{cd}$.  If $J$ intersects $C_c$,
we continue to use $J_c$ for $J\cap C_c$.  Corollary \ref{unionofH} shows
that $J_c$ is a regular $\J$-class of $C_c$ and that $G_J$ is a maximal
subgroup of $J_c$, that is we may
take $G_J = G_{J_c}$.

Let $s,t\in C_c$. Suppose $J$ is a regular $\J$-class of $C^{cd}$
that doesn't intersect $C_c$. Clearly if $x,y\in J$ are idempotents,
then $xsy\notin J$ and $xty\notin J$. Thus
\begin{equation}\label{nointersect}
s\equiv_{(J,G_J,\Rad H {G_J})}t\ \text{if}\ J\cap C_c =\emptyset
\end{equation}
The $\J$-classes of $0$ and $1$, of course, don't intersect $C_c$.

\begin{Lemma}\label{technical}
Let $J$ be a regular $\J$-class of $C^{cd}$ intersecting $C$.
Then
$$\equiv_{(J,G_J,\Rad H {G_J})}{}\cap (C_c\times C_c) =
{}\equiv_{(J_c,G_{J_c},\Rad H
{G_{J_c}})}.$$
\end{Lemma}
\begin{proof}
Choose a Rees matrix representation $\MM^0(G_J,A,B,C)$ for $J^0$.
Since, by Corollary \ref{unionofH},  $J_c$ is a union of $\H$-classes
of $J$, there are subsets $A_c\subseteq A$ and $B_c\subseteq B$
such that $\MM^0(G_J,A_c,B_c,C_c)$ is a Rees matrix representation of
$J_c^0$, where $C_c$ is the restriction of $C$ to $B_c\times A_c$.
Let us denote by
\mbox{$\eta:F(J)\cup J\to \MM^0(G_J,A,B,C)$} the projection.  Clearly the
restriction to \mbox{$F(J_c)\cup J_c$} is the projection
\mbox{$\eta_c:F(J_c)\cup
J_c\to \MM^0(G_J,A_c,B_c,C_c)$}.  Similarly, if
\mbox{$\psi:\MM^0(G_J,A,B,C)\to \MM^0(G_J/\Rad H {G_J},A,B,\ov C)$} is the
natural projection, then its restriction to $\MM^0(G_J,A_c,B_c,C_c)$
is the natural projection \mbox{$\psi_c:\MM^0(G_J,A_c,B_c,C_c)\to
\MM^0(G_J/\Rad H {G_J},A_c,B_c,\ov C_c)$}.

Let $s,t\in C_c$.  Suppose $x,y\in J$ are idempotents.   Since
$x,y\notin 0,1$, they must each belong to some local monoid of $C$.
Suppose at least one of these local monoids is not $C_c$. Then
$xsy=0=xty$ and so
\begin{equation}\label{letsalmostfinish}
xsy\eta\psi = 0
=xty\eta\psi\ \text{if $x$ or $y$ is not in $C_c$}.
\end{equation}
If $x,y\in C_c$, then $x,y\in J_c$.   So $xsy,xty\in F(J_c)\cup J_c$.
Then the previous paragraph shows
\begin{equation}\label{letsfinish}
xsy\eta\psi = xty\eta\psi \iff xsy\eta_c\psi_c = xty\eta_c\psi_c.
\end{equation}
Putting together \eqref{letsalmostfinish}, \eqref{letsfinish} and
\eqref{radcong} (for both $J$ and $J_c$), plus the observation that in
\eqref{radcong} we may
assume that $x,y$ are idempotents, we obtain  $$s \equiv_{(J,G_J,\Rad
  H {G_J})}{}t \iff s\equiv_{(J_c,G_{J_c},\Rad H
{G_{J_c}})}t,$$ thereby completing the proof of the lemma.
\end{proof}

We may now deduce Theorem \ref{supertech}.

\subsection{Proof of Theorem \ref{supertech}}
Suppose that $C$ admits an $\pv {LH}$-morphism $\p:C\to D$ with
$D\in \ell \pv V$.  Then, for $c\in \pvr {Obj}(C)$,
$\p|_{C_c}\to D_c$ is an $\pv {LH}$-morphism by Proposition
\ref{forgotten}.  It follows $C_c\in
\pv {LH}\malce \pv V$ (since $D_c\in \pv V$ by assumption on $D$).

Conversely, suppose that $C\in \ell(\pv {LH}\malce \pv V)$.  The
inclusion map on arrows induces a faithful semigroupoid morphism
$\p_0:C\to C^{cd}$ (it doesn't send local identities to the
identity). Let $\p_1:C^{cd}\rightarrow \prod \GGM (J,G_J,\Rad H {G_J})$ be the
product of the canonical morphisms, where the product runs over all
regular $\J$-classes of $C^{cd}$. Let $\equiv$ be the congruence
on $C$ associated to $\p_0\p_1$, let $D=C/\!{\equiv}$ and let
$\psi:C\to D$ be the quotient map.

\begin{Lemma}\label{finishproof}
The natural projection $\pi:C\to D$ is an $\pv {LH}$-morphism and
$D\in \ell \pv V$.
\end{Lemma}
\begin{proof}
First of all the local monoids of $D$ are of the form $C_c/(\psi)$
with \mbox{$c\in \pvr{Obj}(C)$}.  But the restriction of $\p_0\p_1$ to
$C_c$  is $\p_1|_{C_c}$. However, \eqref{nointersect} and Lemma
\ref{technical} show that the congruence $\sim$ on $C_c$
associated to $\p_1|_{C_c}$ is the intersection of the congruences
$\equiv_{(J_c,G_{J_c},\Rad H {G_{J_c}})}$ associated to the
regular $\J$-classes $J_c$ of $C_c$.  By Theorem \ref{Thm:GGMLH}
and the assumption that $C_c\in \pv {LH}\malce \pv V$, we may
conclude that each $C_c/\!{\equiv_{(J_c,G_{J_c},\Rad H {G_{J_c}})}}\in
\pv V$ and hence $D_c=C_c/\!{\sim}\in \pv V$.  This shows that
$D\in \ell \pv V$.  Also, by Theorem \ref{Prop:LHmap},
$\p_1|_{C_c}$ is an $\pv {LH}$-morphism.  Hence, by Propostion
\ref{forgotten}, $\psi:C\to D$ is an \pv {LH}-morphism, as
desired.
\end{proof}

Lemma \ref{finishproof} completes the proof of Theorem
\ref{supertech} and hence establishes Theorem \ref{mainthm}.\qed

\section*{Acknowledgments}
I would like to thank Bret Tilson for the profound influence he has
had on my work and on my writing of mathematics.  This paper could
not exist if he had not developed his theory of categories as
algebra.

\bibliographystyle{amsplain}

\end{document}